\documentclass[a4paper,12pt]{article}
\usepackage{amssymb,amsmath}

\newcommand{\R}{\mathbb{R}}
\newcommand{\C}{\mathbb{C}}

\newcommand{\ep}{\varepsilon}
\newtheorem{theorem}{Theorem}
\newtheorem{lemma}{Lemma}

\begin{document}
\title{Convexity bounds for $L$-functions}
\author{D.R. Heath-Brown\\Mathematical Institute, Oxford}
\date{}
\maketitle

Let the Dirichlet series
\[F(s)=\sum_{n=1}^{\infty}a_n n^{-s}\]
be absolutely convergent for $\sigma>1$, and extend to a meromorphic
function on $\C$.  Suppose further that $F(s)$ is regular apart
possibly for a pole of order $m\ge 0$ at $s=1$ and that $(s-1)^mF(s)$
is then of finite order.  We assume finally that
$\Phi(s)=\gamma(s)F(s)$ satisfies functional equation
\[\Phi(s)=\overline{\Phi(1-\overline{s})}, \]
where
\[\gamma(s)=\eta Q^s \prod_{j=1}^k\Gamma(\lambda_j s+\mu_j+i\nu_j).\]
Here $\eta\in\C,\, Q\in\R$ and $\lambda_j,\mu_j,\nu_j\in\R\,$ (for $1\le
j\le k$) are constants, satisfying
\[|\eta|=1,\; Q>0,\; \lambda_j>0,\;\mu_j>0.\]
These hypotheses are amongst those required for the ``Selberg Class''
(Selberg \cite{Selb}).

One can now use the Phragm\'{e}n--Lindel\"{o}f theorem in a standard
way to estimate $F(1/2)$ (or, more generally, $F(1/2+it)$).  For
example, if one has $|a_n|\le A(\ep)n^{\ep}$ for any $\ep>0$, then
\begin{equation}\label{conv}
F(1/2)\ll \ep^{-1}A(\ep/2)C^{1/2+\ep},
\end{equation}
where the conductor $C$ is given by
\begin{equation}\label{cond}
C=Q\prod_{j=1}^k(1+|\nu_j|)^{\lambda_j}.
\end{equation}
Here the implied constant depends on $m,k$ and the $\lambda_j$ and
$\mu_j$, but not on $\ep, Q$ or the $\nu_j$.  The aim of this note is
to show how one can remove the extraneous $\ep$ from the exponent 
$1/2+\ep$, under suitable additional hypotheses.  Where one has
appropriate information on the coefficients $a_n$ this can be
done by using some form of approximate functional equation, which will
require an estimate for a sum of the type $\sum_{n\le N}a_n
n^{-1/2}$.  However it is unclear in general how one can bound such
sums efficiently.

Our principal result is the following.
\begin{theorem}
Suppose, in addition to the hypotheses above, that $F(s)$ has an
absolutely convergent Euler product for $\sigma>1$, so that
\begin{equation}\label{bs}
\log F(s)=\sum_n b_n n^{-s}
\end{equation}
with the coefficients $b_n$ supported on the prime powers.  Then
\[F(1/2)\ll C^{1/2}\exp\{4\sum_n|b_n|n^{-3/2}\},\]
with the implied constant depending on $m,k$ and the $\lambda_j$ and
$\mu_j$, but not on $Q$ or the $\nu_j$.
\end{theorem}

The condition that $F(s)$ should have an Euler product is part of the
definition of the Selberg Class.  However all that we require of the
coefficients $b_n$ is that (\ref{bs}) should be absolutely convergent
for $\sigma>1$.  If one were to suppose in addition that
$|b_n|\le cn^{1/3}$, say, then one would of course have a clean bound
$F(1/2)\ll_c C^{1/2}$.

Although our theorem refers only to $F(1/2)$ it is easily modified to
estimate $F(1/2+it)$ in general.  Indeed, if $F(s)$ is entire it
suffices to apply the theorem to $F_t(s)=F(s+it)$, for which
$F_t(1/2)=F(1/2+it)$.  One readily checks here that $F_t(s)$ satisfies
a functional equation of the same form as before, but with the values
$\nu_j$ shifted by $t$.

The proof of our theorem makes it clear that the convexity estimate
above could only be tight if all small non-trivial zeros of $F(s)$
were close to the edge of the critical strip.  (The terms $J(\rho)$
below are genuinely positive for zeros in the interior of the strip.)
\bigskip

A result of the type above, but with stronger hypotheses, was 
described by Soundararajan at the Canadian Number Theory Association 
meeting in Waterloo, and the present paper is an outgrowth of
discussions started there.  Under suitable circumstances
Soundararajan's approach leads \cite{sound1} to an estimate in which 
one saves a further factor of nearly $\log C$.  A weak subconvexity
result of this type is already sufficient for certain applications,
see Holowinsky and Soundararajan \cite{hol}.
The author is grateful to Professor
Soundararajan for a number of interesting remarks, and in particular
for the reference \cite{PS}.
\bigskip

To prove our theorem we will use the following lemma, which can be viewed as a
variant of Jensen's Formula, modified by a conformal transformation so
as to apply to a function in a strip.
\begin{lemma}
Let $P$ be the path
from $\pi/2-i\infty$ to
$\pi/2+i\infty$ and then from $-\pi/2+i\infty$ to $-\pi/2-i\infty$.
When $|\Im(\rho)|<\pi/2$ define
\[J(a+ib)=\log|{\rm coth}(\rho/2)|,\]
and set $J(\rho)=0$ for $|\Im(\rho)|\geq \pi/2$.

Let $f(z)$ be an entire function of finite order, non-vanishing at
$z=0$.  Then
\begin{equation}\label{lem}
\frac{1}{2\pi i}\int_P\log|f(z)|\frac{dz}{\sin z}=\log|f(0)|+
\sum_{\rho}J(\rho),
\end{equation}
where $\rho$ runs over zeros of $f(z)$, counted according to
multiplicity. 
\end{lemma}
Doubtless the hypotheses of this lemma can be weakened
considerably. A related (but different) result is given by 
P\'{o}lya and Szeg\"{o} \cite[pages 119 \& 120]{PS}.

We observe that $J(a+ib)=J(-a+ib)$, and that $J(\rho)\ge 0$ for all
$\rho$. It follows that
\[\log|f(0)|\leq\frac{1}{2\pi}\int_{-\infty}^{\infty}
\log|f(\pi/2+it)f(-\pi/2-it)|\frac{dt}{\cosh t}.\]
Thus after a simple change of variable we find that
if $G$ is entire of finite order then
\[\log|G(1/2)|\leq\frac{1}{2\pi}\int_{-\infty}^{\infty}
\log|G(1+\delta+i\kappa t)G(-\delta+i\kappa t)|\frac{dt}{\cosh t}\]
where we have set $\kappa=(1+2\delta)/\pi$ for convenience.  
We shall apply this with
$G(s)=F(s)(s-1)^m$ and $0<\delta<1$.  The contribution on the right
hand side from terms involving $\log|(s-1)^m|$ is $O(m)$.  Applying
the functional equation, and noting that
\[\int_{-\infty}^{\infty}\frac{dt}{\cosh t}=\pi,\]
leads to a bound
\begin{eqnarray*}
\log|F(1/2)|&\leq&\frac{1}{\pi}\int_{-\infty}^{\infty}
\log|F(1+\delta+i\kappa t)|\frac{dt}{\cosh t}+
(\frac{1}{2}+\delta)\log Q\\
&&\hspace{1cm}\mbox{}+\sum_{j=1}^k
\frac{1}{2\pi}\int_{-\infty}^{\infty}
\log\left|\frac{\Gamma(\alpha_j+i\nu_j(t))}{\Gamma(\beta_j+i\nu_j(t))}\right|
\frac{dt}{\cosh t}+O(m),
\end{eqnarray*}
where we have written
\[\alpha_j=\lambda_j(1+\delta)+\mu_j,\;\;\;
\beta_j=-\lambda_j\delta+\mu_j,\;\;\;
\nu_j(t)=\nu_j+\kappa t\]
for convenience.  However
\[\log\left|\frac{\Gamma(\alpha+i\nu)}{\Gamma(\beta+i\nu)}\right|
=(\alpha-\beta)\log(1+|\nu|)+O_{\alpha,\beta}(1)\]
for $\alpha,\beta>0$, and
\[\int_{-\infty}^{\infty}\log(1+|\nu+\kappa t|)\frac{dt}{\cosh t}
=\pi\log(1+|\nu|)+O(1).\]
We therefore conclude that
\[\log|F(1/2)|\leq\frac{1}{\pi}\int_{-\infty}^{\infty}
\log|F(1+\delta+i\kappa t)|\frac{dt}{\cosh t}+
(\frac{1}{2}+\delta)\log C+O(1),\]
where $C$ is given by (\ref{cond}), and
the implied constant depends on $m,k$ and the $\lambda_j$ and
$\mu_j$, but not on $\delta, Q$ or the $\nu_j$.

It is now easy to deduce (\ref{conv}), but to do better we shall need
the Euler product.  Since
\[\int_{-\infty}^{\infty}n^{-i\kappa t}\frac{dt}{\cosh t}=
\frac{2\pi}{n^{\pi\kappa/2}+n^{-\pi\kappa/2}}=
\frac{2\pi}{n^{1/2+\delta}+n^{-1/2-\delta}}\]
we find that
\[\log|F(1/2)|\leq 2\sum_n \Re(b_n)(n^{3/2+2\delta}+n^{1/2})^{-1}+
(\frac{1}{2}+\delta)\log C+O(1).\]
The theorem now follows on allowing $\delta$ to tend downwards to zero.
\bigskip

We end by giving a sketch proof of the lemma.  
Consider the integral
\[\frac{1}{2\pi i}\int_P\log|1-z/\rho|\frac{dz}{\sin z}.\]
This vanishes if $|\Im(\rho)|\ge\pi/2$, and otherwise is $J(\rho)$.
To see this one sets $\rho=a+ib$ and
integrates around the path $P$ supplemented by a
loop integral from $\pi/2+ib$ around $\rho$ and back. The function
$J(\rho)$ then arises as
\[J(\rho)=2\cosh b\int_a^{\pi/2}\frac{\sin x}{\cosh(2b)-\cos(2x)}dx.\]

The lemma now
follows in the case in which $f(z)$ is a polynomial.  
Moreover if $g(z)$ is also a polynomial then
\[\frac{1}{2\pi i}\int_P\log|\exp(g(z))|\frac{dz}{\sin z}=\Re(g(0)).\]

Now, since $f(z)$ is a function of finite order there is a
positive integer $M$ such that the sum $\sum |\rho|^{-M}$ is 
convergent, and such that
\[f(z)=\exp(h(z))\prod_{\rho}E_M(z/\rho),\]
where $h(z)$ is a suitable polynomial and
\[E_N(w)=(1-w)\exp\{\sum_{j=1}^M w^j/j\}.\]
We may then write $f(z)=f_1(z;N)f_2(z;N)$ with
\[f_1(z;N)=\exp(g(z))\prod_{|\rho|\le N}E_M(z/\rho),\;\;\;
f_2(z;N)=\prod_{|\rho|>N}E_M(z/\rho).\]
We have already shown that the lemma holds for the function $f_1(z;N)$.
Moreover, a crude estimate, considering separately the cases $|z|\le|\rho|/2$,
$|\rho|/2<|z|<2|\rho$ and $|z|\ge 2|\rho|$, shows that
\[\int_P\log|E_M(z/\rho)|\frac{dz}{\sin z}\ll_M |\rho|^{-M}.\]
Thus the contribution to (\ref{lem}) corresponding to 
$\log|f_2(z;N)|$ tends to zero as $N$ goes
to infinity, and the lemma follows.

\bigskip
\bigskip

Mathematical Institute,

24--29, St. Giles',

Oxford

OX1 3LB

UK
\bigskip

{\tt rhb@maths.ox.ac.uk}

\end{document}